\documentclass{amsart}
\usepackage{amsmath}
\usepackage{amsthm}
\usepackage{amssymb}
\usepackage{amsfonts,mathrsfs}
\usepackage{stmaryrd}
\usepackage{amsxtra}  
\usepackage{epsfig}
\usepackage{verbatim}

\theoremstyle{plain}
\newtheorem{theorem}[equation]{Theorem}

\newtheorem{lemma}[equation]{Lemma}
\newtheorem{corollary}[equation]{Corollary}

\theoremstyle{remark}

\numberwithin{equation}{section}

\newcommand{\ssubset}{\subset\subset}
\newcommand{\dbar}{\bar \partial}

\begin{document}

\bibliographystyle{plain}
\title[Plurisubharmonic defining functions]{A note on plurisubharmonic defining functions in $\mathbb{C}^{2}$}
\author{J. E. Forn\ae ss, A.-K. Herbig}
\thanks{Research of the first author was partially supported by an NSF grant}
\address{Department of Mathematics, \newline University of Michigan, Ann Arbor, Michigan 48109}
\email{fornaess@umich.edu}
\address{Department of Mathematics, \newline University of Michigan, Ann Arbor, Michigan 48109}
\email{herbig@umich.edu}
\date{}
\begin{abstract}
  Let $\Omega\ssubset\mathbb{C}^{2}$ be a smoothly bounded domain. Suppose that $\Omega$ admits 
  a  smooth defining function which is plurisubharmonic on the boundary of $\Omega$. Then the 
  Diederich-Forn\ae ss exponent can be chosen arbitrarily close to $1$, and the closure of $\Omega$ 
  admits a Stein neighborhood basis.
\end{abstract}
\maketitle

\section{Introduction}
  Let $\Omega\ssubset\mathbb{C}^{2}$ be a smoothly bounded domain. Throughout, we suppose that $
  \Omega$ admits a $\mathcal{C}^\infty$-smooth defining function $\rho$ which is  plurisubharmonic on 
  the 
  boundary, $b\Omega$, 
  of $\Omega$, 
  i.e.,
   \begin{align*}
      H_{\rho}(\xi,\xi)(z):=
       \sum_{j,k=1}^{2}\frac{\partial^{2}\rho}{\partial z_{j}\partial\bar{z}_{k}}(z)\xi_{j}\bar{\xi}_{k}\geq 0
       \end{align*}
   for all $z\in b\Omega$ and $\xi=(\xi_{1},\xi_{2})\in\mathbb{C}^{2}$. This property  comes 
   up naturally as a sufficiency condition for global regularity of the Bergman projection, see
   \cite{BoaStr91,HerMcN06}.

 The purpose of this paper is to investigate how the plurisubharmonicity of $\rho$ influences the 
  behaviour of the complex Hessian of $\rho$ (or of the complex Hessians of some other defining 
  functions of $\Omega$) away from the boundary of $\Omega$.  
  
  Suppose $D=\{z\in\mathbb{C}^{2}\;|\;r(z)<0\;\}$ is a smoothly bounded, pseudoconvex domain. Then it 
  follows by standard arguments, that 
  there exists a neighborhood $W$ of the boundary of $D$ such 
  that the following lower estimate for the complex Hessian of $r$ holds:
  \begin{align}\label{E:StandardEst}
    H_{r}(\xi,\xi)(q)\geq\mathcal{O}(r(q))|\xi|^{2}
    +\mathcal{O}\left(|\xi|\cdot\left|\langle\partial r(q),\xi\rangle\right|\right)
  \end{align}
  for $q\in W$ and $\xi\in\mathbb{C}^{2}$(see for instance \cite{R1981} for 
  details).
  Our main result shows how to improve the estimate \eqref{E:StandardEst} under the additional 
  condition that there is some smooth defining function of $D$ which is plurisubharmonic on the 
  boundary of $D$.
  
  \begin{theorem}\label{T:MainTheorem}
    Let $\Omega\ssubset\mathbb{C}^{2}$ be a smoothly bounded domain, and suppose $\Omega$ 
    admits a smooth defining function which is plurisubharmonic on the boundary, $b\Omega$, of $ 
    \Omega$.
    Then the following holds:
      for each $\epsilon>0$ and  $K>0$, there exist  a neighborhood $V$ 
      of $b\Omega$  and defining functions $r_{1}$ and $r_{2}$ such that for all 
      $\xi\in\mathbb{C}^{2}$
      \begin{align}\label{E:MainEst}
        H_{r_{1}}(\xi,\xi)(q)\geq\epsilon r_{1}(q)|\xi|^{2}+K|\langle\partial r_{1}(q),\xi\rangle|^{2}
        \;\;\;\text{for}\;\;\; q\in V\cap\overline{\Omega}
      \end{align}
     and
    \begin{align}\label{E:MainEst2}
      H_{r_{2}}(\xi,\xi)(q)\geq -\epsilon r_{2}(q)|\xi|^{2}+K|\langle\partial r_{2},\xi\rangle|^{2}
      \;\;\;\text{for}\;\;\; q\in V\cap\overline{\Omega^{C}}.
    \end{align}
  \end{theorem}
  
   An immediate consequence of Theorem \ref{T:MainTheorem} is the existence of strictly 
   plurisubharmonic exhaustion functions of $\Omega$ and of the complement of $\overline{\Omega}$:
  
   \begin{corollary}\label{C:DFexponent}
    Suppose the hypotheses of Theorem \ref{T:MainTheorem} holds. Then
    \begin{enumerate}
    \item[(i)]  for any $\eta\in(0,1)$ there exists a smooth defining function $\widetilde{r}_{1}$ such that 
       $-(-\widetilde{r}_{1})^{\eta}$ is strictly plurisubharmonic on $\Omega$,
    \item[(ii)] for any $\eta>1$ there exist a neighborhood $V$ of $b\Omega$ and a smooth defining   
       function $\widetilde{r}_{2}$ such that 
       $\widetilde{r}_{2}^{\eta}$ is strictly plurisubharmonic on $V\setminus \overline{\Omega}$.
     \end{enumerate}
  \end{corollary}
  
  A Diederich-Forn\ae ss exponent of a domain $D\ssubset\mathbb{C}^{n}$ is a number 
  $\tau\in(0,1]$ for which there exists a smooth defining function $s$ of $D$ so that $-(-s)^{\tau}$ is strictly 
  plurisubharmonic on $D$. That all 
  smoothly bounded pseudoconvex domains in $\mathbb{C}^{n}$ have a 
  Diederich-Forn\ae ss exponent $\tau$ was shown in \cite{DF1977A} (see also \cite{R1981}). It is also 
  known that there are pseudoconvex domains for which the largest possible $\tau$ 
  might be arbitrarily close to 0 (see \cite{DF1977B}).
  However, part (i) of Corollary \ref{C:DFexponent} says that 
  the Diederich-Forn\ae ss exponent 
  can be chosen arbitrarily close to 1 on domains which admit a smooth 
  defining function, which is plurisubharmonic on the boundary.
  Part (ii) of Corollary \ref{C:DFexponent} is of interest, since it implies that the closure of 
  $\Omega$ has a Stein neighborhood basis. In particular, it follows that $b\Omega$ is uniformly 
  H-convex. We remark that partial results regarding the existence of a Stein neighborhood basis for  
  the closure of  a domain, which satisfies the hypotheses of Theorem \ref{T:MainTheorem}, have been 
  obtained
   in \cite{Sah06}.
  
  \medskip

  The paper is structured as  follows. In Section 2, we identify the obstruction to \eqref{E:MainEst} to  
  hold for the given defining function $\rho$. We then give an example to show 
  that this obstruction might actually occur.
  In Section \ref{S:Modification}, we prove Theorem \ref{T:MainTheorem}, and we conclude this paper by 
  proving Corollary \ref{C:DFexponent} in Section  
  \ref{S:DFexponent} .  
  
  We would like to thank J.D. McNeal for stimulating discussions about this project.
 
  \section{The obstruction}\label{S:Obstruction}
  Throughout, $(z_{1},z_{2})$ will denote the coordinates of $\mathbb{C}^2$. We shall identify
   the vector $\langle\xi_1,\xi_2\rangle$ in $\mathbb{C}^{2}$ with $\xi_1\frac{\partial}{\partial z_1}+
   \xi_2\frac{\partial}{\partial z_2}$ in the $(1,0)$-tangent bundle of $\mathbb{C}^2$ at any given point.
   We use the
   pointwise hermitian inner product $\langle .,.\rangle$  defined by
   $\langle\frac{\partial}{\partial z_j},\frac{\partial}{\partial z_k}\rangle=\delta^j_k$.  We also shall use 
   $\langle.,.\rangle$ to denote contractions of 
   vector fields and forms. We hope this abuse of notation will not confuse the reader as it should be clear 
   from the context what is meant.
  
  \medskip
  
  Let us first see which quantities the right hand side of \eqref{E:StandardEst} depends on. To do so, 
  we need to use Taylor's formula: 
  
  \subsection{Taylor's formula in our context}\label{SS:Taylor}
  Since $b\Omega$ is smooth, there exist a neighborhood $U$ of $b\Omega$ and a smooth map
  \begin{align*}
  \pi:\overline{\Omega}\cap U&\longrightarrow b\Omega\\
  q&\longmapsto \pi(q)=p
  \end{align*}
   such that $p\in b\Omega$ lies on the line normal to $b\Omega$ passing through 
  $q$, and $|p-q|$ is equal to  the complex euclidean distance, $d_{b\Omega}(q)$, of $q$ to 
   $b\Omega$.
  Denote by $\vec{n}_{p}$ the unit outward normal to $b\Omega$ at $p$. Then 
  $q=p-d_{b\Omega}(q)\vec{n}_{p}$. Note that in complex 
  notation
  \begin{align*}
    \vec{n}_{p}=\frac{\left\langle \frac{\partial\rho}{\partial\overline{z}_{1}}, 
    \frac{\partial\rho}{\partial\overline{z}_{2}}\right\rangle}
    {|\partial\rho|}(p),\;\;\text{which implies}\;\;
    q=p-\frac{d_{b\Omega}(q)}{|\partial\rho|}\left\langle\frac{\partial\rho}{\partial \overline{z}_{1}},
    \frac{\partial\rho}{\partial\overline{z}_{2}}\right\rangle(p).
  \end{align*}
  Let $f\in C^{2}(\overline{\Omega})$, $q\in\overline{\Omega}\cap U$ and $p=\pi(q)$. Then Taylor's formula 
  in complex notation says
  \begin{align*}
    f(q)&=f(p)+\sum_{j=1}^{2}\left[ \frac{\partial f}{\partial z_{j}}(p)(q_{j}-p_{j})
    +\frac{\partial f}{\partial\overline{z}_{j}}(p)(\overline{q}_{j}-\overline{p}_{j})\right]
    +\mathcal{O}(|q-p|^{2})\\
    &=
    f(p)-\frac{d_{b\Omega}(q)}{|\partial\rho(p)|}
    \sum_{j=1}^{2}\left[\frac{\partial\rho}{\partial\overline{z}_{j}}(p)\frac{\partial f}{\partial z_{j}}(p)
    +\frac{\partial\rho}{\partial z_{j}}(p)\frac{\partial f}{\partial\overline{z}_{j}}(p)
    \right]+\mathcal{O}(d^{2}_{b\Omega}(q)). 
  \end{align*}
  Define the vector field $N(z)=\frac{1}{|\partial\rho(z)|}\sum_{j=1}^{2}
  \frac{\partial\rho}{\partial\overline{z}_{j}}(z)
  \frac{\partial}{\partial z_{j}}$. Then
  \begin{align}\label{E:Taylor}
    f(q)=f(p)-2d_{b\Omega}(q)\left[(ReN)(f)\right](p)+\mathcal{O}(d_{b\Omega}^{2}(q)).
  \end{align}
  
  \medskip
  \subsection{Partial Taylor analysis of the complex Hessian of $\rho$}\label{SS:Tayloronrho}
  
  After possibly shrinking the neighborhood $U$ of $b\Omega$,  the smooth vector fields
  \begin{align*}
    L=\frac{\frac{\partial\rho}{\partial z_{2}}\frac{\partial}{\partial z_{1}}
    -\frac{\partial\rho}{\partial z_{1}}\frac{\partial}{\partial z_{2}}}
    {|\partial\rho|}\;\;\text{and}\;\;
    N=\frac{\frac{\partial\rho}{\partial\overline{z}_{1}}\frac{\partial}{\partial z_{1}}
    +\frac{\partial\rho}{\partial\overline{z}_{2}}\frac{\partial}{\partial z_{2}}}
    {|\partial\rho|}
  \end{align*}
  are defined on $\overline{\Omega}\cap U$, and it holds that
  \begin{align*}
    L(\rho)=\langle L,N\rangle=0\;\;\text{and}\;\;|L|=1=|N|\;\;\text{on}\;\;\overline{\Omega}\cap U.
  \end{align*}
  Before we get down to business, we need some more notation:
  for vector fields $X(z)=\sum_{i=1}^{2}X_{i}(z)\frac{\partial}{\partial z_{i}}$  and
  $Y(z)=\sum_{i=1}^{2}Y_{i}(z)\frac{\partial}{\partial z_{i}}$, we shall write
  \begin{align*}
     H_{\rho}(X,Y)(z)&=\sum_{j,k=1}^{2}\frac{\partial^{2}\rho}{\partial z_{j}\partial\bar{z}_{k}}(z)
     X_{j}(z)\overline{Y}_{k}(z).
 \end{align*}
 We denote by $\Omega_{W}$ the set of all points $q\in\Omega\cap U$ for which $p=\pi(q)$ is a weakly  
 pseudoconvex boundary point.
 
 Let $\epsilon>0$ be fixed. For each fixed  $q\in\Omega_{W}\cap U$  and $\xi\in\mathbb{C}^{2}$ there  
 exist constants $a_{q,\xi}$ and $b_{q,\xi}$ such that
 $\xi=a_{q,\xi}L(q)+b_{q,\xi}N(q)$. Note that then $|\xi|^{2}=|a_{q,\xi}|^{2}+|b_{q,\xi}|^{2}$. For now, we  
 only consider $q\in\Omega_{W}\cap U$, and for notational ease, we shall drop the subscripts $q,\xi$. 
 We first note that
 \begin{align}\label{E:xiLN}
   H_{\rho}(\xi,\xi)(q)=|a|^{2}H_{\rho}(L,L)(q)+2Re\left(a\overline{b}H_{\rho}(L,N)(q)\right)
   +|b|^{2}H_{\rho}(N,N)(q).
 \end{align}
 We apply \eqref{E:Taylor} to $H_{\rho}(L,L)(q)$, i.e.,
 \begin{align*}
   H_{\rho}(L,L)(q)
   =H_{\rho}(L,L)(p)-2d_{b\Omega}(q)\left(ReN\right)\left(H_{\rho}(L,L)\right)(p)
   +\mathcal{O}(d^{2}_{b\Omega}(q)).
 \end{align*}
 Since $H_{\rho}(L,L)$ is real valued and $H_{\rho}(L,L)(p)=0$, it follows that
 \begin{align*}
   H_{\rho}(L,L)(q)=-2d_{b\Omega}(q)Re(NH_{\rho}(L,L))(p)+\mathcal{O}(d^{2}_{b\Omega}(q)).
 \end{align*}
 Notice that $NH_{\rho}(L,L)(p)$ is real. The last equation combined with \eqref{E:xiLN} gives us then
 \begin{align*}
   H_{\rho}(\xi,\xi)(q)\geq&
   |a|^{2}\left[-2d_{b\Omega}(q)\left(NH_{\rho}(L,L)\right)(p)+\mathcal{O}(d^{2}_{b\Omega}(q))\right]\\
   &-2|a||b||H_{\rho}(L,N)(q)|+|b|^{2}H_{\rho}(N,N)(q).
 \end{align*}
 The Cauchy-Schwarz inequality implies 
 \begin{align*}
   H_{\rho}(\xi,\xi)(q)\geq&
   |a|^{2}\left[-2d_{b\Omega}(q)\left(NH_{\rho}(L,L)\right)(p)-\rho^{2}(q)
   +\mathcal{O}(d^{2}_{b\Omega}(q))\right]\\
   &+|b|^{2}\left[\frac{-1}{\rho^{2}(q)}|H_{\rho}(L,N)(q)|^{2}+H_{\rho}(N,N)(q)\right].
 \end{align*}
 Notice that, after  possibly shrinking  the neighborhood $U$ of $b\Omega$, we can assume that
 \begin{align*}
   -\rho^{2}(q)+\mathcal{O}(d^{2}_{b\Omega}(q))\geq\frac{\epsilon}{4}\rho(q)
 \end{align*}
 for all $q\in\Omega_{W}\cap U$. Therefore,
 \begin{align*}
  H_{\rho}(\xi,\xi)(q)\geq&
   |a|^{2}\left[-2d_{b\Omega}(q)\left(NH_{\rho}(L,L)\right)(p)+\frac{\epsilon}{4}\rho(q)\right]\\
   &+|b|^{2}\left[\frac{-1}{\rho^{2}(q)}|H_{\rho}(L,N)(q)|^{2}+H_{\rho}(N,N)(q)\right]
 \end{align*}
 for all $q\in\Omega_{W}\cap U$.
 Because of the plurisubharmonicity of $\rho$ on $\overline{\Omega}_{W}\cap b\Omega$, it follows that
 \begin{align*}
   |H_{\rho}(L,N)|\leq\left(H_{\rho}(L,L)\right)^{\frac{1}{2}}\left(H_{\rho}(N,N)\right)^{\frac{1}{2}}
 \end{align*}
 holds on $\overline{\Omega}_{W}\cap b\Omega$. Since $q\in\Omega_{W}\cap U$, i.e., since $\pi(q)=p$ is a weakly pseudoconvex boundary point, we get that $H_{\rho}(L,N)(p)=0$. Therefore, there exists a constant $c_{1}>0$, depending on $\rho$, such that
 \begin{align*}
   |H_{\rho}(L,N)(q)|^{2}\leq c_{1}|\rho(q)|^{2}\;\;\text{for all}\;\;q\in\Omega_{W}\cap U.
 \end{align*}
 This gives us the following lower bound on $H_{\rho}(\xi,\xi)(q)$:
 \begin{align*}
   H_{\rho}(\xi,\xi)(q)\geq& |a|^{2}\left[-2d_{b\Omega}(q)\left(NH_{\rho}(L,L)\right)(p)
   +\frac{\epsilon}{4}\rho(q)\right]\\
   &-|b|^{2}\left[c_{1}-H_{\rho}(N,N)(q)\right],
  \end{align*}
  which implies that for some constant $c_{2}>0$ depending on $\rho$
  \begin{align}\label{E:TayloronHrho} 
   H_{\rho}(\xi,\xi)(q)\geq |a|^{2}\left[-2d_{b\Omega}(q)\left(NH_{\rho}(L,L)\right)(p)
   +\frac{\epsilon}{4}\rho(q)\right]-c_{2}|b|^{2}
 \end{align}
 holds for $q\in\Omega_{W}\cap U$.
 
 Note that \eqref{E:TayloronHrho} is a more detailed version of \eqref{E:StandardEst} for those points 
 $q\in\Omega$ near $b\Omega$ whose projections $\pi(q)$ are weakly pseudoconvex boundary points.
 Moreover, 
 inequality \eqref{E:MainEst} is within range, if $NH_{\rho}(L,L)$ is non-positive at all weakly 
 pseudconvex boundary points. The term  $NH_{\rho}(L,L)$ being positive at some weakly pseudoconvex boundary point $p_{0}$ means that the function $H_{\rho}(L,L)$ decreases when one moves from $p_{0}$ inside the domain along the line normal to $b\Omega$ at $p_{0}$. This, of course, means that $H_{\rho}(L,L)$ becomes negative there, which destroys any hope of $\rho$ being plurisubharmonic in some neighborhood of $p_{0}$. Clearly, $NH_{\rho}(L,L)(p_{0})> 0$ obstructs
 inequality \eqref{E:MainEst} to hold for all $\epsilon>0$.
 \medskip
 
 \subsection{Example \& idea of modification of $\rho$} We shall first give an example of a domain where $NH_{\rho}(L,L)$ is positive at a  weakly pseudoconvex boundary point. Consider the domain $D=\{(z,w)\in\mathbb{C}^{2}\;|\;\rho(z,w)<0\}$ near the origin, where
 \begin{align*}
   \rho(z,w)=Re(w)+|w|^{2}+Re(w)|z|^{2}+|z|^{2}|w|^{2}+|z|^{4}+|z|^{6}.
 \end{align*}
 One can easily show that $\rho$ is plurisubharmonic on $bD$ near the origin. In fact, $\rho$ is strictly plurisubharmonic on $bD$ near the origin except when $z=0$. Let $\xi=(\xi_{1},\xi_{2})$ and $q=(0,w)$ be a point in $D$ which lies on the line normal to $bD$ through the origin. Then
 \begin{align*}
 H_{\rho}(\xi,\xi)(q)=(Re(w)+|w|^{2})|\xi_{1}|^{2}+|\xi_{2}|^{2}.
 \end{align*}
 Thus $\rho$ can not be plurisubharmonic in any neighborhood of the origin. Note that 
 this is caused by the term $Re(w)|z|^{2}$ contained in the definition of $\rho$. However, this is our old enemy, that is
 \begin{align*}
   (NH_{\rho}(L,L))(0)=\frac{\partial}{\partial w}\left(
   \frac{\partial^{2}\rho}{\partial z\partial\bar{z}}
   \right)(0)=\frac{1}{2}>0!
 \end{align*}
 Now the question is, whether we can manipulate $\rho$ such that the obstruction vanishes. Notice that 
 the answer to that in the above example is yes: let $r(z,w)=\rho(z,w)/(1+|z|^{2})$. Then
 \begin{align*}
   r(z,w)=Re(w)+|w|^{2}+|z|^{4},
 \end{align*}
 which is plurisubharmonic everywhere.
 \medskip
 
 Recall, that we actually want to show an estimate like
 \begin{align*}
  -2d_{b\Omega}(q)(NH_{\rho}(L,L))(p)\geq\epsilon\rho(q).
 \end{align*}
 Obviously, just multiplying $\rho$ by a small positive number is not going to remove the obstruction.
 So, we consider another defining function $\rho\cdot h$ of $\Omega$, where $h$ is some smooth, 
 positive function. We shall now list a few characteristics of $h$ which should give us some control on 
 the obstruction term:
 \begin{enumerate}
   \item In order to use the basic estimate \eqref{E:TayloronHrho} for $\rho\cdot h$, we would need that $ 
     \rho\cdot h$ is still plurisubharmonic at weakly pseudoconvex boundary points. This can be achieved, 
     if we choose $h$ such that all its first order derivatives vanish at all weakly pseudoconvex boundary 
     points.
   \item We need to consider those third order derivatives of $\rho\cdot h$, which are forced upon us by 
   the
   obstruction term, at weakly pseudoconvex points. If we assume that all first order derivatives of $h$ 
   vanish at weakly pseudoconvex 
   points (and if we ignore, at least temporarily, that in the obstruction term $N$ does not only act on the 
   Levi form of $\rho\cdot h$ but  also on $L$ and $\overline{L}$), then there are only two terms to be
   considered:
   \begin{enumerate}
     \item There is the product of the original obstruction term of $\rho$ and $h$, which tells us that $h$ 
     itself should not be large at the weakly pseudoconvex points.
     \item There are the terms which involve one derivative of $\rho$ and two derivatives of $h$. Since we 
     are on $b\Omega$ the only such term which can appear is $N$ acting on $\rho$ multiplied 
     with the Levi form of $h$. This seems to say that we need the Levi form of $h$ to be negative definite 
     at the weakly pseudoconvex points. One can show that
     $NH_{\rho}(L,L)$ equals $LH_{\rho}(N,L)$ at weakly pseudoconvex points 
     (see \eqref{E:something1} and the following lemma).  Thus the obstruction term itself gives us a 
     function, $-|H_{\rho}(N,L)|^{2}$, whose Levi form is strictly negative definite at those points where 
     it is needed.
   \end{enumerate}
 \end{enumerate}
 Clearly, we can not choose $h$ to be $-|H_{\rho}(N,L)|^{2}$, since the latter function vanishes at weakly pseudoconvex points, and hence $\rho\cdot h$ would not be a defining function of $\Omega$. However, taking (1) and (2) into account $e^{-|H_{\rho}(N,L)|^{2}}$ seems like a suitable candidate for $h$.
  
 \section{Proof of Theorem \ref{T:MainTheorem}}\label{S:Modification}
 Let $C>0$ be a large constant, which will be chosen later. We will consider the smooth defining 
 function
 \begin{align*}
   r_{C}=r=\rho e^{-C\sigma},\;\;\text{where}\;\;\sigma=|H_{\rho}(N,L)|^{2},
 \end{align*}
 and we shall work with the vector fields
 \begin{align*}
  L^{r}=\frac{\frac{\partial r}{\partial z_{2}}\frac{\partial}{\partial z_{1}}
    -\frac{\partial r}{\partial z_{1}}\frac{\partial}{\partial z_{2}}}
    {|\partial r|}\;\;\text{and}\;\;
    N^{r}=\frac{\frac{\partial r}{\partial\overline{z}_{1}}\frac{\partial}{\partial z_{1}}
    +\frac{\partial r}{\partial\overline{z}_{2}}\frac{\partial}{\partial z_{2}}}
    {|\partial r|},
 \end{align*}
 which are defined on $\overline{\Omega}\cap U$ (after possibly shrinking $U$). As before, we note that $L^{r}(r)=\langle L^{r},N^{r}\rangle=0$ and $|L^{r}|=|N^{r}|=1$. Moreover, on $b\Omega$ we have 
 $L^{r}=L$ and $N^{r}=N$. 
 
 As before, we suppose that $q\in\Omega_{W}\cap U$. Here, the decomposition of a vector
 $\xi\in\mathbb{C}^{2}$ with respect to the vector fields $L^{r}$ and $N^{r}$ at $q$ is different than 
 before. Clearly, we can write $\xi=a_{q,\xi}L^{r}(q)+b_{q,\xi}N^{r}(q)$ again; however, the constants 
 $a_{q,\xi}$ and $b_{q,\xi}$ are different from before. Again, for notational convenience, we shall drop 
 those subscripts $q,\xi$. 
 
 Let us first see whether the basic estimate \eqref{E:TayloronHrho} holds for $r$.
 The only special property of $\rho$, which we used to derive \eqref{E:TayloronHrho}, is that 
 $H_{\rho}(L,N)(p)=0$, where $p$ is a weakly pseudoconvex boundary point. Thus, to see whether
 \eqref{E:TayloronHrho} holds for $r$ we shall compute $H_{r}(L^{r},N^{r})(p)$. A straightforward computation gives
 \begin{align*}
 \frac{\partial ^{2}r}{\partial z_{j}\partial\overline{z}_{k}}
 =
 e^{-C\sigma}\left[-C\frac{\partial\sigma}{\partial\overline{z}_{k}}
 \left(\frac{\partial \rho}{\partial z_{j}}
 -C\rho\frac{\partial\sigma}{\partial z_{j}}\right)\right.
 &+\frac{\partial^{2}\rho}{\partial z_{j}\partial\overline{z}_{k}}\\
 &\left.-C\frac{\partial\rho}{\partial\overline{z}_{k}}
 \frac{\partial\sigma}{\partial z_{j}}-C\rho\frac{\partial^{2}\sigma}{\partial z_{j}\partial\overline{z}_{k}}
 \right].
 \end{align*}
 Since $H_{\rho}(L,N)(p)=0$, it follows that not only $\sigma$ but also any derivative of $\sigma$ at $p$ 
 vanishes, and thus we obtain
 \begin{align*}
    \frac{\partial ^{2}r}{\partial z_{j}\partial\overline{z}_{k}}(p)= 
    \frac{\partial ^{2}\rho}{\partial z_{j}\partial\overline{z}_{k}}(p).
 \end{align*}
 In particular, $r$ is plurisubharmonic at $p$ and $H_{r}(L^{r},N^{r})(p)=0$. Thus \eqref{E:TayloronHrho}
 holds for $r$. That is: there exists a constant $c_{2}>0$ (depending on $r$) such that
  \begin{align}\label{E:TayloronHr}
    H_{r}(\xi,\xi)(q)\geq |a|^{2}\left[-2d_{b\Omega}(q) \left(N^{r}H_{r}(L^{r},L^{r})\right)(p)+
     \frac{\epsilon}{4}r(q)\right]-c_{2}|b|^{2}
   \end{align}  
   holds for all $q\in\Omega_{W}\cap U$ after possibly shrinking $U$.
   \medskip
   
   To see whether we truly gain anything by using $r$ instead of $\rho$, we have to figure out how
   $(N^{r}H_{r}(L^{r},L^{r}))(p)$ is related to $(NH_{\rho}(L,L))(p)$. We shall prove the following
   \begin{align}\label{Claim}
     \text{\underline{Claim}:}\;\;
     N^{r}H_{r}(L^{r},L^{r})(p)\leq \left[NH_{\rho}(L,L)-C|\partial\rho|\cdot(NH_{\rho}(L,L))^{2}\right](p).
   \end{align}
  Note that  $N^{r}=N$ on $b\Omega$, which implies on $b\Omega$
 \begin{align*}
   N^{r}H_{r}(L^{r},L^{r})
   =
   NH_{r}(L^{r},L^{r})
   =
   \sum_{\ell=1}^{2}N_{l}\frac{\partial}{\partial z_{\ell}}
   \left(
   \sum_{j,k=1}^{2}\frac{\partial ^{2}r}{\partial z_{j}\partial\overline{z}_{k}}L^{r}_{j}
   \overline{L}^{r}_{k}
   \right).
 \end{align*}
 Since $L^{r}$ is a weak complex tangential direction at $p$ and $r$ is plurisubharmonic at $p$, 
 we have 
 \begin{align*}
   \sum_{j,k=1}^{2}\frac{\partial^{2} r}{\partial z_{j}\overline{z}_{k}}
   \left(\sum_{\ell=1}^{2}N_{\ell}\frac{\partial L^{r}_{j}}{\partial z_{\ell}}\right)\overline{L}^{r}_{k}(p)
   =0=
   \sum_{j,k=1}^{2}\frac{\partial^{2} r}{\partial z_{j}\overline{z}_{k}} L_{j}^{r}
   \left(\sum_{\ell=1}^{2}N_{\ell}\frac{\partial \overline{L}^{r}_{k}}{\partial z_{\ell}}\right)(p).
 \end{align*}
 Moreover, we have  $L^{r}(p)=L(p)$, which gives us that
 \begin{align*}
   \left(N^{r}H_{r}(L^{r},L^{r})\right)(p)
   =\left(\sum_{j,k,\ell=1}^{2}\frac{\partial^{3}r}{\partial z_{j}\partial\overline{z}_{k}\partial z_{\ell}}
   L_{j}\overline{L}_{k}N_{\ell}\right)(p).
 \end{align*} 
Let us now compute those third derivatives of $r$: 
 \begin{align*}
 &\frac{\partial^{3} r}{\partial z_{j}\partial\overline{z}_{k}\partial z_{\ell}}\\
 =&
 e^{-C\sigma}\left[-C\frac{\partial\sigma}{\partial z_{\ell}}
 \left\{
 -C\frac{\partial\sigma}{\partial\overline{z}_{k}}
 \left(
 \frac{\partial\rho}{\partial z_{j}}-C\rho\frac{\partial\sigma}{\partial z_{j}}
 \right)
 +\frac{\partial^{2}\rho}{\partial z_{j}\partial\overline{z}_{k}}
 \right.
 -C\frac{\partial\rho}{\partial\overline{z}_{k}}\frac{\partial\sigma}{\partial z_{j}}
 -C\rho\frac{\partial^{2}\sigma}{\partial z_{j}\overline{z}_{k}}\right\}\\
 &\hspace{1cm}-C\frac{\partial^{2}\sigma}{\partial\overline{z}_{k}z_{\ell}}
 \left(\frac{\partial\rho}{\partial z_{j}}-C\rho\frac{\partial\sigma}{\partial z_{j}}\right)
 -C\frac{\partial\sigma}{\partial\overline{z}_{k}}
 \left(\frac{\partial^{2}\rho}{\partial z_{j}\partial z_{\ell}}
 -C\frac{\partial\rho}{\partial z_{\ell}}\frac{\partial\sigma}{\partial z_{j}}-C\rho
 \frac{\partial^{2}\sigma}{\partial z_{j}\partial z_{\ell}}
 \right)\\
 &\hspace{0.8cm}\left.+\frac{\partial^{3}\rho}{\partial z_{j}\partial\overline{z}_{k}\partial z_{\ell}}
 -C\left(\frac{\partial^{2}\rho}{\partial\overline{z}_{k}\partial z_{\ell}}
 \frac{\partial\sigma}{\partial z_{j}}+\frac{\partial\rho}{\partial\overline{z}_{k}}
 \frac{\partial^{2}\sigma}{\partial z_{j}\partial z_{\ell}}
 +\frac{\partial\rho}{\partial z_{\ell}}\frac{\partial^{2}\sigma}{\partial z_{j}\overline{z}_{k}}
 +\rho\frac{\partial^{3}\sigma}{\partial z_{j}\partial\overline{z}_{k}\partial z_{\ell}}
 \right)\right].
 \end{align*}
 First note that $\rho$ as well as $\sigma$ and all its first order derivatives vanish at $p$. Also, since $L$ is complex tangential to $b\Omega$ all the terms involving $\frac{\partial\rho}{\partial z_{j}}$ or
$\frac{\partial\rho}{\partial\overline{z}_{k}}$ vanish. Thus we get 
  \begin{align*}
  \left(N^{r}H_{r}(L^{r},L^{r})\right)(p)
  =
  \left(NH_{\rho}(L,L)-C\langle\partial\rho,N\rangle H_{\sigma}(L,L)
  \right)(p).
 \end{align*}
 Since $\langle\partial\rho, N\rangle(p)=|\partial\rho(p)|$, it follows that
 \begin{align*}
   \left(N^{r}H_{r}(L^{r},L^{r})\right)(p)
  =
  \left(NH_{\rho}(L,L)-C|\partial\rho|H_{\sigma}(L,L)
  \right)(p).
 \end{align*}
 Recall that $\sigma=|H_{\rho}(N,L)|^{2}$.  Using that $H_{\rho}(N,L)(p)=0$, a direct computation gives 
 us
 \begin{align*}
   H_{\sigma}(L,L)(p)&=
   \left|\langle \partial H_{\rho}(N,L),L\rangle(p)\right|^{2}
   +
   \left|\langle \overline{\partial} H_{\rho}(N,L),\overline{L}\rangle(p)\right|^{2}\\
   &\geq
   \left|\langle \partial H_{\rho}(N,L),L\rangle(p)\right|^{2}.
 \end{align*}
 We compute further
 \begin{align*}
   \langle\partial H_{\rho}(N,L),L\rangle
   &=
   \sum_{j=1}^{2}L_{j}\frac{\partial}{\partial z_{j}}
   \left(
   \sum_{k,\ell=1}^{2}\frac{\partial^{2}\rho}{\partial z_{\ell}\overline{z}_{k}}
   N_{\ell}\overline{L}_{k}
   \right)\\
   &=
   \sum_{j,k,\ell=1}^{2}
   \frac{\partial^{3}\rho}{\partial z_{j}\partial\overline{z}_{k}\partial z_{\ell}}
   L_{j}\overline{L}_{k}N_{\ell}
   +\sum_{k,\ell=1}^{2}\frac{\partial^{2}\rho}{\partial z_{l}\partial\overline{z}_{k}}
   \left(\sum_{j=1}^{2}L_{j}\frac{\partial}{\partial z_{j}}
   \left(\overline{L}_{k}N_{\ell}
   \right)
   \right).
 \end{align*}
 Since $L$ is a weak complex tangential direction at $p$ and $\rho$ is plurisubharmonic at $p$, it follows that
 \begin{align}\label{E:something1}
 \langle\partial H_{\rho}(N,L),L\rangle(p)
 =
 NH_{\rho}(L,L)(p)+
 \sum_{k,\ell=1}^{2}
 \frac{\partial^{2}\rho}{\partial z_{\ell}\partial\overline{z}_{k}}N_{\ell}
 \left(
 \sum_{j=1}^{2}L_{j}\frac{\partial \overline{L}_{k}}{\partial z_{j}}
 \right)(p).
 \end{align}
We claim that the last term on the right hand side vanishes:
\begin{lemma}\label{L:Zisweak}
  Suppose $X$ is a smooth vectorfield, which is complex tangential to $b\Omega$. 
  Furthermore, suppose that  $b\Omega$ is weakly pseudoconvex at some boundary point $p$. Define
  $Y=\sum_{j=1}^{2}\overline{X}_{j}\frac{\partial X_{k}}{\partial \overline{z}_{j}}
  \frac{\partial}{\partial z_{k}}$.
  Then $Y$ is weak complex tangential to $b\Omega$ at $p$.
\end{lemma}
\begin{proof}
  Since $X$ is tangential, $X(\rho)=0$ holds on $b\Omega$. Moreover, we 
  have $\overline{X}(X(\rho))=0$ on $b\Omega$. Therefore
  \begin{align*}
    0&=\overline{X}(X(\rho))(p)
    =\sum_{j,k=1}^{2}\overline{X}_{j}\frac{\partial}{\partial\overline{z}_{j}}
    \left(
    X_{k}\frac{\partial \rho}{\partial z_{k}}
    \right)(p)\\
    &=\sum_{j,k=1}^{2}\overline{X}_{j}\frac{\partial X_{k}}{\partial \overline{z}_{j}}
  \frac{\partial \rho}{\partial z_{k}}(p)
  +
  \sum_{j,k=1}^{2}\frac{\partial^{2}\rho}{\partial z_{k}\partial\overline{z}_{j}}
  X_{k}\overline{X}_{j}(p)=Y(\rho)(p),
  \end{align*}
  where the last step holds since $H_{\rho}(X,X)(p)=0$ by our hypothesis. Thus, $Y$ is complex 
  tangential direction at $p$. In particular, $H_\rho(Y,Y)(p)=0$.
 \end{proof}
 If we set $X=L$, Lemma \ref{L:Zisweak} implies that the last term in
 \eqref{E:something1} vanishes.
 Thus, we obtain
 \begin{align*}
 H_{\sigma}(L,L)(p)
 \geq 
 \left|
 \langle
 \partial H_{\rho}(N,L),L
 \rangle(p)
 \right|^{2}=
 \left|
 NH_{\rho}(L,L)(p)
 \right|^{2},
 \end{align*}
 which proves the Claim \eqref{Claim}. That is
 \begin{align*}
   N^{r}H_{r}(L^{r},L^{r})(p)\leq\left[NH_{\rho}(L,L)-C|\partial\rho|\cdot(NH_{\rho}(L,L))^{2}\right](p).
 \end{align*}
 Hence, the lower estimate \eqref{E:TayloronHr} on the complex Hessian of $r$ now becomes
 \begin{align}
    H_{r}(\xi,\xi)(q)
   \geq&
   |a|^{2}
   \left[
   2d_{b\Omega}(q)\left\{Cc_{3}\left(NH_{\rho}(L,L)\right)^{2}
   -NH_{\rho}(L,L)\right\}(p)+
   \frac{\epsilon}{4}r(q)
   \right]\notag\\
   &-c_{2}|b|^{2}
 \end{align}
 for  $q\in\Omega_{W}\cap U$, where $c_{3}>0$ is chosen such that $|\partial\rho|\geq c_{3}$ on $b\Omega$.

 \medskip
 
 We are now set to show that there exist a $C>0$ and a neighborhood $U_{C}$ of $b\Omega$ such that
 \begin{align}\label{E:toshow}
 2d_{b\Omega}(q)\left[
 C c_{3}\left( NH_{\rho}(L,L)\right)^{2}
 - NH_{\rho}(L,L)
 \right](p)\geq\frac{\epsilon}{4}r(q)
 \end{align}
 holds for $q\in\Omega_{W}\cap U_{C}$, which would imply that \eqref{E:MainEst} holds for these points.
 
  To make our life easier, let us write
 $A_{p}$ for $NH_{\rho}(L,L)(p)$, i.e., \eqref{E:toshow} becomes
 \begin{align*}
   2d_{b\Omega}(q)\left[C c_{3}A_{p}^{2}-A_{p}\right]\geq\frac{\epsilon}{4}r(q).
 \end{align*}
 If $C c_{3}A_{p}^{2}-A_{p}\geq 0$, then \eqref{E:toshow} holds trivially. Moreover, increasing $C$ does not destroy this non-negativity.
 Suppose that $C c_{3}A_{p}^{2}-A_{p}<0$. First notice that there exists a constant $c_{4}>0$ such that
  $d_{b\Omega}(q)\leq c_{4}|\rho(q)|$
 for all $q\in\Omega\cap U$. Since $\rho=re^{C\sigma}$, it follows that
  $d_{b\Omega}(q)\leq c_{4}e^{C\sigma(q)}|r(q)|$.
 Thus, to prove \eqref{E:toshow} it is sufficient to show
 \begin{align*}
   2c_{4}e^{C\sigma(q)}|r(q)|\left[C c_{3}A_{p}^{2}-A_{p}\right]&\geq \frac{\epsilon}{4}r(q),\;\;
   \text{which is equivalent to}\\
   e^{C\sigma(q)}\left[C c_{3}A_{p}^{2}-A_{p}\right]&\geq-\frac{\epsilon}{8c_{4}}.
 \end{align*}
 Let $U_{C}\subset U$ be a neighborhood of $b\Omega$ such that 
 $z\in\Omega\cap U_{C}$ implies that $e^{C\sigma(z)}\leq 2e^{C\sigma(\pi(z))}$. Notice that $U_{C}$ is a true neighborhood of $b\Omega$, since $\sigma$ is smooth near $b\Omega$. Moreover, in the situation which we are considering, i.e., where $\pi(q)$ is a weakly pseudoconvex boundary point, we then have that $q\in\Omega_{W}\cap
 U_{C}$ implies $e^{C\sigma(q)}\leq 2$.
 Therefore, to obtain \eqref{E:toshow} it is sufficient that
 \begin{align*}
   Cc_{3}A_{p}^{2}-A_{p}\geq-\frac{\epsilon}{16 c_{4}}
 \end{align*}
 holds on $\Omega_{W}\cap U_{C}$. We remark that neither $c_{3}, c_{4}$ nor $A_{p}$ depend on the 
 choice of  $C$. Thus, choosing
  \begin{align*}
    C=\max\left\{0,\max_{p\in b\Omega\;\text{weak}}\frac{\frac{-\epsilon}{16c_{4}}+A_{p}}{c_{3}A_{p}^{2}} 
    \right\}
  \end{align*}
  proves \eqref{E:toshow} on $\Omega\cap U_{C}$, which implies that
  \begin{align}\label{E:Estonweak}
    H_{r}(\xi,\xi)(q)\geq \frac{\epsilon}{2}r(q)|\xi|^{2}-c_{2}|\langle\partial r(q),\xi\rangle|^{2}
  \end{align}
  holds on $\Omega_{W}\cap U_{C}$.
  
  \medskip
  
  Let us show now that an estimate similar to \eqref{E:Estonweak} holds  near $\Omega_W\cap U_C$. 
  Note first that our computations leading up to
  \eqref{E:Estonweak} imply that $N^r H_r(L^r,L^r)\leq\frac{\epsilon}{16}$ holds on the set of the weakly 
  pseudoconvex boundary points of $\Omega$. Hence by continuity, there exists a neighborhood
  $W$ of the set of weakly pseudoconvex boundary points such that 
  $N^rH_r(L^r,L^r)\leq\frac{\epsilon}{8}$ on $W\cap b\Omega$. We may assume
  that $W\subset U_C$ and that $q\in W\cap\Omega$ implies $\pi(q)\in W\cap b\Omega$.
  Using Taylor's formula, it follows 
  for $q\in W\cap\Omega$ that
  \begin{align*}
    H_r(L^r,L^r)(q)&\geq
    H_r(L^r,L^r)(\pi(q))+
    \frac{\epsilon}{4}r(q)+\mathcal{O}(r^2(q))\\
    &\geq H_r(L^r,L^r)(\pi(q))+\frac{\epsilon}{2}r(q)
  \end{align*}
  after possibly shrinking of $W$. 
  Another application of Taylor's formula gives us for $q\in W\cap\Omega$
  \begin{align*}
    H_r(\xi,\xi)(q)\geq&
    |a|^2\left[H_r(L^r,L^r)(\pi(q))+\frac{\epsilon}{2}r(q)\right]
    +
    |b|^2H_r(N^r,N^r)\\
    &\hspace{3.35cm}+2|a||b|\left[|H_r(L^r,N^r)(\pi(q))|+\mathcal{O}(r(q))\right]\\
    \geq&
     |a|^2\left[H_r(L^r,L^r)(\pi(q))+\epsilon r(q)\right]-c_5|\langle\partial r(q),\xi\rangle|^2\\
     &\hspace{3.3cm}+2|a||b||H_r(L^r,N^r)(\pi(q))|,
  \end{align*}
  where the last step follows by the Cauchy-Schwarz inequality for some constant
   $c_5>0$ sufficiently large. Since $H_r(L^r,L^r)(\pi(q))$ is positive, we only need to estimate the term
   $|H_r(L^r,N^r)(\pi(q))|$.
  Note  first that $r$ is not plurisubharmonic on $b\Omega$ at strictly 
  pseudoconvex boundary points, though $\rho$ is. However, since any derivative of $\sigma$
  is $\mathcal{O}(H_\rho(N^r,L^r))$ on $b\Omega$ and since $\rho$ is plurisubharmonic on $b\Omega$,  
  it follows that there exists a constant $c_6>0$ such that  
  \begin{align*}
    |H_r(L^r,N^r)|^2\leq c_6H_r(L^r,L^r)\left[H_r(N^r,N^r)+c_6\right]\;\;\text{on}\;\;b\Omega.
  \end{align*}
  The Cauchy-Schwarz inequality implies now, that for some constant $c_7>0$ we have
  \begin{align}\label{E:EstonnbhdW}
    H_r(\xi,\xi)(q)\geq
    \epsilon r(q)|\xi|^2-c_7\left|\langle\partial r(q),\xi\rangle\right|^2
  \end{align}
  for $q\in W\cap \Omega$.
    We define $r_{1}=r+Kr^{2}$ for some $K>2c_{7}$. Note that
     \begin{align*}
       H_{r_{1}}(\xi,\xi)(q)=(1+2Kr)H_{r}(\xi,\xi)(q)+2K|\langle\partial r,\xi\rangle|^{2}.
     \end{align*}
     Let $U_{K}=\{\;z\in W\;|\;1+2Kr(z)\geq\frac{1}{2}\;\}$, then \eqref{E:EstonnbhdW} implies 
     for $q\in\Omega\cap U_{K}$
     \begin{align*}
       H_{r_{1}}(\xi,\xi)(q)&\geq\frac{1}{2}\epsilon r(q)|\xi|^{2}
       +(2K-c_{7})|\langle\partial r(q),\xi\rangle|^{2}\\
       &\geq \epsilon r_{1}(q)|\xi|^{2}+K|\langle\partial r_{1}(q),\xi\rangle|^{2}.
     \end{align*}
     That is, we have shown that \eqref{E:MainEst} holds on $\Omega\cap U_K$.

     Next we show that \eqref{E:MainEst} also holds near the remaining strictly pseudoconvex boundary 
     points. We note that $S=b\Omega\setminus(W \cap b\Omega)$ is a closed 
     subset of the set of the strictly pseudoconvex boundary points. This implies, as long as $K>0$ is 
     chosen
     sufficiently large,  that there exists a neighborhood $U_{S}$ of $S$ such that $r_{1}$ is strictly 
     plurisubharmonic on $\Omega\cap U_{S}$. In particular, there exists a neighborhood $V$ of 
     $b\Omega$ such that
     \begin{align*}
      H_{r_{1}}(\xi,\xi)(q)\geq\epsilon r_{1}(q)|\xi|^{2}+K|\langle\partial r_{1}(q),\xi\rangle|^{2}
     \end{align*}
     for all $q\in\Omega\cap V$ and $\xi\in\mathbb{C}^{2}$. This proves \eqref{E:MainEst}.
     
     \medskip
     
     The proof of \eqref{E:MainEst2} is very similar to the above proof of \eqref{E:MainEst}. In fact, we only 
     need to change a few signs to derive \eqref{E:MainEst2}. First, one realizes that the basic estimate
     \eqref{E:TayloronHrho} becomes: there exist a neighborhood $U$ of $b\Omega$ and a constant 
     $c_{2}>0$ such that
     \begin{align*}
       H_{\rho}(\xi,\xi)(q)\geq |a|^{2}\left[ 2d_{b\Omega}(q)(NH_{\rho}(L,L)(p)-\frac{\epsilon}{4}\rho(q)\right]
       -c_{2}|b|^{2}
     \end{align*}
     for $q\in\Omega^{C}\cap U$. Here, as before, the points $q$ in consideration are such that their 
     orthogonal projections $\pi(q)=p$
     onto $b\Omega$ are weakly pseudoconvex boundary points. As one would expect, we have an 
     obstruction to plurisubharmonicity of $\rho$ outside of $\overline{\Omega}$ at those weakly 
     pseudoconvex boundary points where $H_{\rho}(L,L)$ decreases along the outward normal. Thus, 
     we should multiply $\rho$ by a smooth, positive function which is strictly plurisubharmonic at those 
     boundary points where $NH_{\rho}(L,L)$ is negative, i.e., we work with the function
     $r=\rho e^{C|H_{\rho}(N,L)|^{2}}$ for $C>0$. Using arguments analog to the ones in the proof of 
     \eqref{E:MainEst}, one can then show that for any $\epsilon>0$ and $K>0$, there exist a 
     neighborhood $V$ of $b\Omega$ and a constant $C>0$ such that the complex Hessian of
     $r_{2}=r+Kr^{2}$ satisfies \eqref{E:MainEst2} on $\Omega^{C}\cap V$.

     \medskip
     \section{Proof of Corollary \ref{C:DFexponent}}\label{S:DFexponent}
     In the following section, we give the proof of Corollary \ref{C:DFexponent}. We start out with part (i) by 
     showing  first that for 
     any $\eta\in(0,1)$ there exist a $\delta>0$, a smooth defining function $r$ of $\Omega$ and a 
     neighborhood $W$ of $b\Omega$  such that
     $g_{1}=-(-re^{-\delta|z|^{2}})^{\eta}$ is strictly plurisubharmonic on $\Omega\cap W$.
     
     Let $\eta\in(0,1)$ be fixed, and $r$ be a smooth defining function of $\Omega$. For notational ease 
     we write $\phi=\delta|z|^{2}$ for $\delta>0$. Here, $r$ and $\delta$ are fixed and to be chosen later.
     Let us compute the complex Hessian of $g_{1}$ on $\Omega\cap W$:
     \begin{align*}
       H_{g_{1}}(\xi,\xi)=&\eta(-r)^{\eta-2}e^{-\phi\eta}
       \left[
       (1-\eta)\right.\left|\langle\partial r,\xi\rangle\right|^{2}-rH_{r}(\xi,\xi)\\
       &+2r\eta Re\left(\langle\partial r,\xi\rangle\langle\overline{\partial\phi,\xi}\rangle\right)
       \left.
       -r^{2}\eta\left|\langle\partial\phi,\xi\rangle\right|^{2}
      +r^{2}H_{\phi}(\xi,\xi)\right].
     \end{align*}
     An application of the Cauchy-Schwarz inequality gives
     \begin{align*}
       2r\eta Re\left(\langle\partial r,\xi\rangle\langle\overline{\partial\phi,\xi}\rangle\right)
       \geq -(1-\eta)\left|\langle\partial r,\xi\rangle\right|^{2}
       -\frac{r^{2}\eta^{2}}{1-\eta}\left|\langle\partial\phi,\xi\rangle\right|^{2}.
     \end{align*}
     Therefore, we obtain for the complex Hessian of $g$ on $\Omega$ the following:
     \begin{align*}
       H_{g_{1}}(\xi,\xi)
       \geq
       \eta(-g_{1})(-r)^{-1}\left[H_{r}(\xi,\xi)
       +
       (-r)\left\{H_{\phi}(\xi,\xi)-\frac{\eta}{1-\eta}\left|\langle\partial\phi,\xi\rangle\right|^{2}
       \right\}
       \right].
     \end{align*}
     Notice that
     \begin{align*}
       H_{\phi}(\xi,\xi)-\frac{\eta}{1-\eta}\left|\langle\partial\phi,\xi\rangle\right|^{2}
       &=
       \delta\left(H_{|z|^{2}}(\xi,\xi)-\frac{\eta}{1-\eta}\delta\left|\langle\overline{z},
       \xi\rangle\right|^{2}
       \right)\\
       &\geq
       \delta|\xi|^{2}\left(
       1-\frac{\eta D}{1-\eta}\delta
       \right),
     \end{align*}
     where $D:=\max_{z\in\overline{\Omega}}|z|^{2}$. Now set 
     $\delta=\frac{1-\eta}{2\eta D}$; it is noteworthy that $\delta$ goes to $0$ as $\eta$ approaches 
     $1^{-}$. We now have
     \begin{align*}
       H_{\phi}(\xi,\xi)-\frac{\eta}{1-\eta}\left|\langle\partial\phi,\xi\rangle\right|^{2}
       \geq\frac{\delta}{2}|\xi|^{2},
     \end{align*}
     which implies that
     \begin{align}\label{E:generalDFest}
       H_{g_{1}}(\xi,\xi)\geq \eta(-g_{1})(-r)^{-1}
       \left[
       H_{r}(\xi,\xi)+\frac{\delta}{2}(-r)|\xi|^{2}
       \right]
     \end{align}
     holds on $\Omega$. 
     
     By \eqref{E:MainEst} 
     there exist a neighborhood $W$ of $b\Omega$ and a smooth defining function $r_{1}$ 
     of $\Omega$ such that
     \begin{align*}
     H_{r_{1}}(\xi,\xi)(q)
     \geq
     \frac{\delta}{4} r_{1}(q)|\xi|^{2}
     \end{align*}
     for all $q\in\Omega\cap W$. Setting $r=r_{1}$ and using \eqref{E:generalDFest}, we obtain
     \begin{align*}
       H_{g_{1}}(\xi,\xi)(q)\geq
       \eta(-g_{1}(q))\cdot\frac{\delta}{4}|\xi|^{2}\;\;\text{for}\;\;q\in\Omega\cap W.
     \end{align*}
     It follows by standard arguments that there exists a defining function $\widetilde{r}_{1}$ such that 
     $-(-\widetilde{r}_{1})^{\eta}$ is strictly 
     plurisubharmonic on $\Omega$; for details see pg. 133 in \cite{DF1977A}. 
      This proves part (i) of Corollary \ref{C:DFexponent}.
     
     \medskip
     
     The proof of part (ii) is similar to the proof of part (i).  Let 
     $\eta>1$ be fixed. We would like to show  that there exists a neighborhood $V$ of $b\Omega$ such 
     that $g_{2}=(re^{\delta|z|^{2}})^{\eta}$ is strictly 
     plurisubharmonic on
     $\overline{\Omega}^{C}\cap V$ for some smooth defining function $r$ and some constant 
     $\delta>0$. Let $W$ be a neighborhood of $b\Omega$.
     Choose $\delta=\frac{\eta-1}{2\eta D}$, where $D=\max_{z\in\overline{W}}|z|^{2}$. Then calculations 
     similar to the ones in the proof of part (i) yield
     \begin{align*}
       H_{g_{2}}(\xi,\xi)\geq \eta g_{2}r^{-1}\left[
       H_{r}(\xi,\xi)+\frac{\delta}{2}r|\xi|^{2}
       \right]\;\;\text{on}\;\;\overline{\Omega}^{C}\cap W.
     \end{align*}
     By \eqref{E:MainEst2} there exist a neighborhood $V$ of $b\Omega$ and a smooth defining function
     $r_{2}$ of $\Omega$ such that
     \begin{align*}
       H_{r_{2}}(\xi,\xi)(q)\geq-\frac{\delta}{4}r_{2}(q)|\xi|^{2}
     \end{align*}
      for all $q\in\overline{\Omega}^{C}\cap V$. Since we may assume that $V\subset W$, it follows that
      \begin{align*}
        H_{g_{2}}(\xi,\xi)(q)\geq\eta g_{2}(q)\cdot\frac{\delta}{4}|\xi|^{2}\;\;\text{for}\;\;
        q\in\overline{\Omega}^{C}\cap V,
      \end{align*}
      which proves \eqref{E:MainEst2}.

\end{document}